\documentclass[12pt,thmsa]{article}
\usepackage{amsfonts}
\usepackage{amsmath}
\usepackage{amssymb}

\input{tcilatex}
\begin{document}

\title{Abstract Harmonic Analysis and Group Algebra on the Group $AN$ . }
\author{Kahar El-Hussein \\
\textit{Department of Mathematics, Faculty of Science, }\\
\textit{\ Al Furat University, Dear El Zore, Syria \ \ and}\\
\textit{Department of Mathematics, Faculty of Arts Science Al Quryyat, }\\
\textit{\ Al-Jouf University, KSA }\\
\textit{E-mail : kumath@ju.edu.sa, kumath@hotmail.com }}
\maketitle

\begin{abstract}
Let $G$ $=SL(n,\mathbb{R})$ be the real semisimple Lie group and let $SL(n,%
\mathbb{R})=KAN$ be the Iwasawa decomposition. The purpose of this paper is
to define the Fourier transform in order to obtain the Plancherel theorem on
the group $AN.$ Besides, we prove the existence theorems for the invariant
differential operators on the nilpotent group $N,$and the solvable $AN.$ To
this end, we give a classification of all left ideals of the group algebra $%
L^{1}(AN).$
\end{abstract}



\section{\textbf{\ Introduction.}}

\bigskip\ \textbf{1.1. }Abstract harmonic analysis is one of the more modern
branches of harmonic analysis, having its roots in the mid-twentieth
century, is analysis on topological groups. The core motivating idea are the
various Fourier transforms, which can be generalized to a transform of
functions defined on Hausdorff locally compact topological groups. If the
group is neither abelian nor compact, no general satisfactory theory is
currently known. In mathematics and theoretical physics, the invariant
differential operators play very important role in the description of
physical symmetries. starting from the early occurrences in the Maxwell,
d'Allembert, Dirac, equations. Usually, the action of the group has the
meaning of a change of coordinates (change of observer) and the invariance
means that the operator has the same expression in all admissible
coordinates. For example the Dirac operator in physics is invariant with
respect to the Poincar\'{e} group and Maxwell equation is invariant by the
Lorentz group .

\textbf{1.2.} Thus, it is important for the applications in mathematical
physics to study systematically such operators and their solutions. Our
interest is to consider the solvable Lie group $AN$ of the real semisimple
Lie group $G=SL(n,\mathbb{R})$, Therefore let $SL(n,\mathbb{R})=KAN$ be the
Iwasawa decomposition. In this paper, we will define the Fourier transform
on the group $AN$ to obtain the Plancherel formula and to solve any
invariant differential operators on the group $AN$. In the end, we give a
classification of all left ideals of group algebra $L^{1}(AN)$ and we
introduce some new groups

\section{\protect\bigskip Notation and Results}

\textbf{2.1. }The fine structure of the nilpotent Lie groups will help us to
do the Fourier transform on a simply connected nilpotent Lie groups $N.$ As
well known any group connected and simply connected $N$ has the following
form%
\begin{equation}
N=\left( 
\begin{array}{cccccccccccc}
1 & x_{1}^{1} & x_{1}^{2} & x_{1}^{3} & . & . & . & . & . & x_{1}^{n-2} & 
x_{1}^{n-1} & x_{1}^{n} \\ 
0 & 1 & x_{2}^{2} & x_{2}^{3} & . & . & . & . & . & x_{2}^{n-2} & x_{2}^{n-1}
& x_{2}^{n} \\ 
0 & 0 & 1 & x_{3}^{3} & . & . & . & . & . & x_{3}^{n-2} & x_{3}^{n-1} & 
x_{3}^{n} \\ 
0 & 0 & 0 & 1 & . & . & . & . & . & x_{4}^{n-2} & x_{4}^{n-1} & x_{4}^{n} \\ 
. & . & . & . & . & . & . & . & . & . & . & . \\ 
. & . & . & . & . & . & . & . & . & . & . & . \\ 
. & . & . & . & . & . & . & . & . & . & . & . \\ 
. & . & . & . & . & . & . & . & . & . & . & . \\ 
. & . & . & . & . & . & . & . & . & x_{n-2}^{n-2}. & .x_{n-2}^{n-1} & 
x_{n-2}^{n} \\ 
0 & 0 & 0 & 0 & . & . & . & . & . & 1 & x_{n-1}^{n-2} & x_{n=1}^{n} \\ 
0 & 0 & 0 & 0 & . & . & . & . & . & 0 & 1 & x_{n}^{n} \\ 
0 & 0 & 0 & 0 & . & . & . & . & . & 0 & 0 & 1%
\end{array}%
\right)
\end{equation}

\bigskip As shown, this matrix is formed by the subgroup $\mathbb{R}$, $%
\mathbb{R}^{2}$,...., $\mathbb{R}^{n-1}$, and $\mathbb{R}^{n}$ 
\begin{equation}
\left( \mathbb{R=}\left[ 
\begin{array}{c}
x_{1}^{1} \\ 
1 \\ 
0 \\ 
0 \\ 
. \\ 
. \\ 
. \\ 
. \\ 
. \\ 
0 \\ 
0 \\ 
0%
\end{array}%
\right] ,\mathbb{R}^{2}=\left[ 
\begin{array}{c}
x_{1}^{2} \\ 
x_{2}^{2} \\ 
1 \\ 
0 \\ 
. \\ 
. \\ 
. \\ 
. \\ 
. \\ 
0 \\ 
0 \\ 
0%
\end{array}%
\right] ,..,\mathbb{R}^{n-1}=\left[ 
\begin{array}{c}
x_{1}^{n-1} \\ 
x_{2}^{n-1} \\ 
x_{3}^{n-1} \\ 
x_{4}^{n-1} \\ 
. \\ 
. \\ 
. \\ 
. \\ 
.x_{n-2}^{n-1} \\ 
x_{n-1}^{n-2} \\ 
1 \\ 
0%
\end{array}%
\right] ,\mathbb{R}^{n}=\left[ 
\begin{array}{c}
x_{1}^{n} \\ 
x_{2}^{n} \\ 
x_{3}^{n} \\ 
x_{4}^{n} \\ 
. \\ 
. \\ 
. \\ 
. \\ 
x_{n-2}^{n} \\ 
x_{n=1}^{n} \\ 
x_{n}^{n} \\ 
1%
\end{array}%
\right] \right)
\end{equation}%
$\ \ \ \ \ \ \ \ \ \ \ \ \ \ \ \ \ \ \ \ $

Each $\mathbb{R}^{i}$ is a subgroup of $N$ of dimension $i$ , $1\leq i\leq
n, $ put $d=n+(n-1)+....+2+1=\frac{n(n+1)}{2},$ which is the dimension of $N$
. According to $[9],$ the group $N$ is isomorphic onto the following group

\begin{equation}
(((((\mathbb{R}^{n}\rtimes _{\rho _{n}})\mathbb{R}^{n-1})\rtimes _{\rho
_{n-1}})\mathbb{R}^{n-2}\rtimes _{\rho _{n-2}}.....)\rtimes _{\rho _{2}}%
\mathbb{R}^{2})\rtimes _{\rho _{1}}\mathbb{R}
\end{equation}

\bigskip That means

\begin{equation}
N\simeq (((((\mathbb{R}^{n}\rtimes _{\rho _{n}})\mathbb{R}^{n-1})\rtimes
_{\rho _{n-1}})\mathbb{R}^{n-2}\rtimes _{\rho _{n-2}}.....)\rtimes _{\rho
_{4}}\mathbb{R}^{3}\rtimes _{\rho _{3}}\mathbb{R}^{2})\rtimes _{\rho _{2}}%
\mathbb{R}
\end{equation}

\textbf{2.2}$.$ Denote by $L^{1}(N)$ the Banach algebra that consists of all
complex valued functions on the group $N$, which are integrable with respect
to the Haar measure of $N$ and multiplication is defined by convolution on $%
N $ as follows:%
\begin{equation}
g\ast f(X)=\int\limits_{N}f(Y^{-1}X)g(Y)dY
\end{equation}%
for any $f\in L^{1}(N)$ and $g\in L^{1}(N),$ where $X=(X^{1},$ $X^{2},$ $%
X^{3},....,X^{n-2},X^{n-1},X^{n}),$ $%
Y=(Y^{1},Y^{2},Y^{3},....,Y^{n-2},Y^{n-1},Y^{n}),$ $X^{1}=x_{1}^{1},$ $%
X^{2}=(x_{1}^{2},x_{2}^{2}),$ $X^{3}=(x_{1}^{3},x_{2}^{3},x_{3}^{3})$ $%
,...., $ $%
X^{n-2}=(x_{1}^{n-2},x_{2}^{n-2},x_{3}^{n-2},x_{4}^{n-2},...,x_{n-2}^{n-2}),$
$%
X^{n-1}=(x_{1}^{n-1},x_{2}^{n-1},x_{3}^{n-1},x_{4}^{n-1},...,x_{n-2}^{n-1},x_{n-1}^{n-1}), 
$ $%
X^{n}=(x_{1}^{n},x_{2}^{n},x_{3}^{n},x_{4}^{n},...,x_{n-2}^{n},x_{n-1}^{n},x_{n}^{n}) 
$ and $dY=dY^{1}dY^{2}dY^{3},....,dY^{n-2}dY^{n-1}dY^{n}$ is the Haar
measure on $N$ and $\ast $ denotes the convolution product on $N.$ We denote
by $L^{2}(N)$ its Hilbert space. Let $\mathcal{U}\;$be the complexified
universal enveloping algebra of the real Lie algebra $\underline{g}$\ of $N$%
; which is canonically isomorphic to the algebra of all distributions on $N$
supported by $\left\{ 0\right\} ,$ where $0$ is the identity element of $N$.
For any $u\in \mathcal{U}$ one can define a differential operator $P_{u}$ on 
$N$ as follows:%
\begin{equation}
P_{u}f(X)=u\ast f(X)=\int\limits_{N}f(Y^{-1}X)u(Y)dY
\end{equation}

\bigskip The mapping $u\rightarrow P_{u}$ is an algebra isomorphism of $%
\mathcal{U}$ onto the algebra of all invariant differential operators on $N$

\section{Fourier Analysis on $N.$}

\bigskip \textbf{3.1. }For each $\ 2\leq i\leq n,$ let 
\begin{eqnarray*}
&&K_{i-1} \\
&=&\left\{ \underbrace{\left\{ \underbrace{\left\{ \underbrace{\left\{ 
\underbrace{\mathbb{R}^{n}}_{\rtimes _{\rho _{n}}}\right\} \times \mathbb{R}%
^{n-1}\times \mathbb{R}^{n-1}}_{\rtimes _{_{\rho _{n-1}}}}\right\} \times 
\mathbb{R}^{n-2}\times \mathbb{R}^{n-2}\times .....}_{\rtimes _{\rho
i+1}}\right\} \times \mathbb{R}^{i}\times \mathbb{R}^{i}}_{\rtimes _{\rho
i}}\right\} \times \mathbb{R}^{i-1}\mathbb{\rtimes }_{\rho _{i}}\mathbb{R}%
^{i-1}
\end{eqnarray*}%
be the group with the following law

\begin{eqnarray}
&&(X....,X^{i},X^{i},U^{i-1},X^{i-1}).(Y....,Y^{i},Y^{i},V^{i-1},Y^{i-1}) \\
&=&[(X....,X^{i},X^{i}).(\rho
_{i}(X^{i})(Y....,Y^{i},Y^{i}))],U^{i-1}+V^{i-1},X^{i-1}+Y^{i-1})  \notag
\end{eqnarray}

For any $X\in L_{i+1},Y\in L_{i+1},$ $X^{i}\in \mathbb{R}^{i},U^{i-1}\in 
\mathbb{R}^{i-1},X^{i-1}\in \mathbb{R}^{i-1},Y^{i}\in \mathbb{R}%
^{i},V^{i-1}\in \mathbb{R}^{i-1},Y^{i-1}\in \mathbb{R}^{i-1}.$ So we get for 
$i=2$

\ 
\begin{eqnarray*}
&&K_{1} \\
&=&\left\{ \underbrace{\left\{ \underbrace{\left\{ \underbrace{\left\{ 
\underbrace{\mathbb{R}^{n}}_{\rtimes _{\rho _{n}}}\right\} \times \mathbb{R}%
^{n-1,1}\rtimes _{\rho _{n-1}}\mathbb{R}^{n-1}}_{\rtimes _{_{\rho
_{n-1}}}}\right\} \times \mathbb{R}^{n-2,1}\rtimes _{\rho _{3}}\mathbb{R}%
^{n-2}\times .....}_{\rtimes _{\rho 3}}\right\} \times \mathbb{R}%
^{2,1}\times \mathbb{R}^{2}}_{\rtimes _{\rho _{2}}}\right\} \times \mathbb{R}%
^{1,1}\mathbb{\rtimes }_{\rho _{2}}\mathbb{R}^{1}
\end{eqnarray*}%
and for $i=n$ , we get 
\begin{equation}
K_{n-1}=\left\{ \underbrace{\mathbb{R}^{n}}_{\rtimes _{\rho _{n}}}\right\}
\times \mathbb{R}^{n-1}\rtimes _{\rho _{n}}\mathbb{R}^{n-1}  \notag
\end{equation}

\bigskip \textbf{3.2. }Let $M=\mathbb{R}^{n}\times \mathbb{R}^{n-1,1}\times 
\mathbb{R}^{n-2,1}\times .....\times \mathbb{R}^{3,1}\times \mathbb{R}%
^{2,1}\times \mathbb{R}^{1,1}\mathbb{=R}^{d}$ be the Lie group, which is the
direct product of $\mathbb{R}^{n,1}=\mathbb{R}^{n},\mathbb{R}^{n-1,1}=%
\mathbb{R}^{n-1},\mathbb{R}^{n-2,1}=\mathbb{R}^{n-2},.....,\mathbb{R}^{3,1}=%
\mathbb{R}^{3},\mathbb{R}^{2,1}=\mathbb{R}^{2}$ and $\mathbb{R}^{1,1}=%
\mathbb{R}^{1}$. Denote by $L^{1}(M)$ the Banach algebra that consists of
all complex valued functions on the group $M$, which are integrable with
respect to the Lebesgue measure on $M$ and multiplication is defined by
convolution on $M$ as: 
\begin{equation}
g\ast _{c}f(X)=\int\limits_{M}f(X-Y)g(Y)dY
\end{equation}%
for any $f\in L^{1}(M),$ $g\in L^{1}(M),$ where $\ast _{c}$ signifies the
convolution product on the abelian group $M.$ We denote again by $\mathcal{U}%
\;$be the complexified universal enveloping algebra of the real Lie algebra $%
\underline{m}$\ of $M$; which is canonically isomorphic to the algebra of
all distributions on $M$ supported by $\left\{ 0\right\} ,$ where $0$ is the
identity element of $M$. For any $u\in \mathcal{U}$ one can define a
differential operator $Q_{u}$ on $M$ as follows:%
\begin{equation}
Q_{u}f(X)=u\ast _{c}f(X)=f\ast _{c}u(X)=\int\limits_{M}f(X-Y)u(Y)dY
\end{equation}

\bigskip The mapping $u\rightarrow Q_{u}$ is an algebra isomorphism of $%
\mathcal{U}$ onto the algebra of all invariant differential operators (the
algebra of all linear differential operators with constant coefficients) on $%
M$

The group $N$ can be identified with the subgroup%
\begin{equation}
N=\left\{ \underbrace{\left\{ \underbrace{\left\{ \underbrace{\left\{ 
\underbrace{\mathbb{R}^{n}}_{\rtimes _{\rho _{n}}}\right\} \times \mathbb{\{}%
0\mathbb{\}}\rtimes _{\rho _{n-1}}\mathbb{R}^{n-1}}_{\rtimes _{_{\rho
_{n-1}}}}\right\} \times \mathbb{\{}0\mathbb{\}}\rtimes _{\rho _{3}}\mathbb{R%
}^{n-2}\times .....}_{\rtimes _{\rho 3}}\right\} \times \mathbb{\{}0\mathbb{%
\}}\times \mathbb{R}^{2}}_{\rtimes _{\rho _{2}}}\right\} \times \mathbb{\{}0%
\mathbb{\}\rtimes }_{\rho _{2}}\mathbb{R}
\end{equation}%
of $K_{1}$ and $M$ can be identified with the subgroup 
\begin{eqnarray*}
&&M \\
&=&\left\{ \underbrace{\left\{ \underbrace{\left\{ \underbrace{\left\{ 
\underbrace{\mathbb{R}^{n}}_{\rtimes _{\rho _{n}}}\right\} \times \mathbb{R}%
^{n-1,1}\mathbb{\rtimes }_{\rho _{n}}\mathbb{\{}0\mathbb{\}}}_{\rtimes
_{_{\rho _{n-1}}}}\right\} \times \mathbb{R}^{n-2,1}\mathbb{\rtimes }_{\rho
_{n-1}}\mathbb{\{}0\mathbb{\}}\times .....}_{\rtimes _{\rho 3}}\right\}
\times \mathbb{R}^{2,1}\mathbb{\rtimes }_{\rho _{3}}\mathbb{\{}0\mathbb{\}}}%
_{\rtimes _{\rho _{2}}}\right\} \times \mathbb{R}^{1,1}\mathbb{\rtimes }%
_{\rho _{2}}\mathbb{\{}0\mathbb{\}}
\end{eqnarray*}

In this paper, we show how the Fourier transform on the vector group $%
\mathbb{R}^{d}$ can be generalized on $N$ and obtain the Plancherel theorem.

\textbf{Definition 3.1.} \textit{For} $1\leq i\leq n,$ \textit{let} $%
\mathcal{F}^{i}$\textit{\ be the Fourier transform on} $\mathbb{R}^{i}$ 
\textit{and }$0\leq j\leq n-1,$ \textit{let }$\dprod\limits_{0\leq l\leq j}%
\mathbb{R}^{n-l}=(..((((\mathbb{R}^{n}\rtimes _{\rho _{n}})\mathbb{R}%
^{n-1})\rtimes _{\rho _{n-1}})\mathbb{R}^{n-2}\rtimes _{\rho
_{n-2}})..\times _{\rho _{n-j}}\mathbb{R}^{n-j}),$ \textit{and let }$%
\dprod\limits_{0\leq l\leq j}\mathcal{F}^{n-l}$ $=\mathcal{F}^{n}\mathcal{F}%
^{n-l}\mathcal{F}^{n-2}....\mathcal{F}^{n-j},$\textit{we can define the
Fourier transform on } $\dprod\limits_{0\leq l\leq n-1}\mathbb{R}^{n-1}=%
\mathbb{R}^{n}\rtimes _{\rho _{n}}\mathbb{R}^{n-1}$ $\rtimes _{\rho _{n-1}}%
\mathbb{R}^{n-2}\rtimes _{\rho _{n-2}}.....\rtimes _{\rho _{3}}\mathbb{R}%
^{2}\rtimes _{\rho _{2}}\mathbb{R}^{1}$\textit{as}%
\begin{eqnarray}
&&\mathcal{F}^{n}\mathcal{F}^{n-1}\mathcal{F}^{n-2}....\mathcal{F}^{2}%
\mathcal{F}^{1}f(\lambda ^{n},\text{ }\lambda ^{n-1},\lambda ^{n-2},\text{%
.......,}\lambda ^{2},\lambda ^{1})  \notag \\
&=&\dint\limits_{N}f(X^{n},X^{n-1},...,X^{2},X^{1})e^{-\text{ }i\langle 
\text{ }(\lambda ^{n},\text{ }\lambda ^{n-1}),(X^{n},X^{n-1}),..,(\lambda
^{2},\text{ }\lambda ^{1}),(X^{2},X^{1})\rangle }\text{ }  \notag \\
&&dX^{n}dX^{n-1}...dX^{2}dX^{1}
\end{eqnarray}%
\textit{for any} $f\in L^{1}(N),$ \textit{where }$%
X=(X^{n},X^{n-1},..,X^{n},X^{n-1}),$ $\mathcal{F}^{d}=\mathcal{F}^{n}%
\mathcal{F}^{n-1}\mathcal{F}^{n-2}....\mathcal{F}^{2}\mathcal{F}^{1}$\ 
\textit{is} \textit{the classical Fourier transform on} $N,$ $%
dX=dX^{n}dX^{n-1}...dX^{2}dX^{1},$ $\lambda =(\lambda ^{n},$ $\lambda
^{n-1},\lambda ^{n-2},..$,$\lambda ^{2},\lambda ^{1}),$ \textit{and}

$\langle (\lambda ^{n},$ $\lambda ^{n-1}),(X^{n},X^{n-1}),..,(\lambda ^{2},$ 
$\lambda ^{1}),(X^{2},X^{1})\rangle =\dsum\limits_{i=1}^{n}X_{i}^{n}\lambda
_{i}^{n}+\dsum\limits_{j=1}^{n-1}X_{j}^{n-1}\lambda
_{j}^{n-1}+..+\dsum\limits_{i=1}^{2}X_{i}^{2}\lambda _{i}^{2}+X^{1}\lambda
^{1}$

\textbf{Definition 3.2. }For every function $f$ on $N$ one can define a
function $\widetilde{f}$ on $K_{1}$ as

\begin{equation}
\widetilde{f}(\rho _{2}(X^{1})[(\rho _{3}(X^{2})(..(\rho
_{n-1}(X^{n-2})((\rho
_{n}(X^{n-1})(X^{n})),U^{i-1}+X^{i-1}),U^{i-2}+X^{i-2})..)),U^{2}+X^{2})],U^{1}+X^{1})
\notag
\end{equation}

\bigskip \textbf{Remark 3.1.} The function $\widetilde{f}$ \ is invariant in
the following sense%
\begin{eqnarray}
&&\widetilde{f}(\rho _{2}(U^{1})[(\rho _{3}(U^{2})(..((\rho
_{n}(U^{n-1})(X^{n})),X^{n-1}-U^{n-1}),X^{n-2}-U^{n-2})..)),X^{2}-U^{2})],X^{1}-U^{1})
\notag \\
&&=\widetilde{f}%
(X,U^{i-1},X^{i-1},U^{i-2},X^{i-2},..,U^{2},X^{2},U^{1},X^{1})
\end{eqnarray}

\bigskip \textbf{Theorem 3.1.} For every function $F$ on $K_{1}$ invariant
in sense and for every function $\phi $ on $N$%
\begin{eqnarray}
&&\phi \ast F(X,U^{n-1},X^{n-1},U^{n-2},X^{n-2},..,U^{2},X^{2},U^{1},X^{1}) 
\notag \\
&=&\phi \text{ }\ast
_{c}F(X,U^{n-1},X^{n-1},U^{n-2},X^{n-2},..,U^{2},X^{2},U^{1},X^{1})
\end{eqnarray}%
\textit{for every} $%
(X,U^{n-1},X^{n-1},U^{n-2},X^{n-2},..,U^{2},X^{2},U^{1},X^{1})$ $\in K_{1}$, 
\textit{where} $\ast $\textit{\ signifies the convolution product on} $N$ 
\textit{with respect the variables} $(X,X^{i-1},X^{i-2},..,X^{2},X^{1})$\ 
\textit{and} $\ast _{c}$\textit{signifies the commutative convolution
product on} $M$ \textit{with respect the variables }$%
(X,,U^{i-1},U^{i-2},..,U^{2},U^{1})$

\textit{Proof: }%
\begin{eqnarray*}
&&\phi \ast F(X,U^{n-1},X^{n-1},U^{n-2},X^{n-2},..,U^{2},X^{2},U^{1},X^{1})
\\
&=&\int\limits_{N}F\left[
(Y^{n},Y^{n-1},Y^{n-2},..,Y^{2},Y^{1})^{-1}(X,U^{n-1},X^{n-1},U^{n-2},X^{n-2},..,U^{2},X^{2},U^{1},X^{1})%
\right] \\
&&\phi (Y^{n},Y^{n-1},Y^{n-2},..,Y^{2},Y^{1})dYdY^{i-1}dY^{i-2}..dY^{2}dY^{1}
\\
&=&\int\limits_{N}F((\rho
(-Y^{1})(Y^{n},Y^{n-1},Y^{n-2},..,Y^{2})^{-1}),-Y^{1})(X,U^{n-1},X^{n-1},U^{n-2},X^{n-2},..,U^{2},X^{2},U^{1},X^{1}))
\\
&&\phi
(Y,Y^{n},Y^{n-1},Y^{n-2},..,Y^{2},Y^{1})dYdY^{n-1}dY^{n-2}..dY^{2}dY^{1} \\
&=&\int\limits_{N}F((\rho
(-Y^{1})((Y^{n},Y^{n-1},Y^{n-2},..,Y^{2})^{-1}(X,U^{n-1},X^{n-1},U^{n-2},X^{n-2},..,U^{2},X^{2},U^{1},X^{1})),X^{1}-Y^{1}))
\\
&&\phi
(Y^{n},Y^{n-1},Y^{n-2},..,Y^{2},Y^{1})dY^{n}dY^{n-1}dY^{n-2}..dY^{2}dY^{1} \\
&=&\int\limits_{N}F(\rho _{2}(-Y^{1})(\rho
_{3}(-Y^{2})((Y^{n},Y^{n-1},Y^{n-2},..,)^{-1}(X^{n},U^{n-1},X^{n-1},U^{n-2},X^{n-2},..,)U^{2},X^{2}-Y^{2},
\\
&&U^{1},X^{1}-Y^{1})]\phi
(Y^{n},Y^{n-1},Y^{n-2},..,Y^{2},Y^{1})dY^{n}dY^{n-1}dY^{n-2}..dY^{2}dY^{1}
\end{eqnarray*}

By a simple calculation, we find

\begin{eqnarray}
&&\phi \ast
F(X^{n},U^{i-1},X^{i-1},U^{i-2},X^{i-2},..,U^{2},X^{2},U^{1},X^{1})  \notag
\\
&=&\int\limits_{N}F\left[
(X^{n}-Y^{n},U^{n-1},X^{n-1}-Y^{n-1},U^{n-2},X^{n-2}-Y^{n-2},..,U^{2},X^{2}-Y^{2},U^{1},X^{1}-Y^{1})%
\right]  \notag \\
&&\phi
(Y^{n},Y^{i-1},Y^{i-2},..,Y^{2},Y^{1})dY^{n}dY^{i-1}dY^{i-2}..dY^{2}dY^{1} 
\notag \\
&=&\int\limits_{M}F\left[
(X^{n}-Y^{n},U^{n-1}-Y^{n-1},X^{n-1},U^{n-2}-Y^{n-2},X^{n-2},..,U^{2}-Y^{2},X^{2},U^{1}-Y^{1},X^{1})%
\right]  \notag \\
&&\phi
(Y^{n},Y^{n-1},Y^{n-2},..,Y^{2},Y^{1})dY^{n}dY^{n-1}dY^{n-2}..dY^{2}dY^{1} 
\notag \\
&=&\phi \ast
_{c}F(X^{n},U^{n-1},X^{n-1},U^{n-2},X^{n-2},..,U^{2},X^{2},U^{1},X^{1})
\end{eqnarray}

\textbf{Corollary. }For every function $F$ $\in $ $L^{1}(K_{1})$ invariant
in sense and for every function $\phi $ $\in $ $L^{1}(N)$%
\begin{eqnarray}
&&\dint \dint ..\dint \dint \mathcal{F}^{d}(\phi \ast F)(\xi ^{n},\lambda
^{n-1},\mu ^{n-1},\lambda ^{n-2},\mu ^{n-2},..,\lambda ^{2},\mu ^{2},\lambda
^{1},\mu ^{1})d\mu ^{n-1}d\mu ^{n-2}..d\mu ^{2}d\mu ^{1}  \notag \\
&=&\mathcal{F}^{d}F(\xi ^{n},\lambda ^{n-1},0,\lambda ^{n-2},0,..,\lambda
^{2},0,\lambda ^{1},0)\mathcal{F}^{d}\phi (\xi ^{n},\lambda ^{n-1},\lambda
^{n-2},..,\lambda ^{2},\lambda ^{1})
\end{eqnarray}

\bigskip \textit{Proof: In fact we get}

\begin{eqnarray}
&&\dint \dint ..\dint \dint \mathcal{F}^{d}(\phi \ast F)(\xi ^{n},\lambda
^{n-1},\mu ^{n-1},\lambda ^{n-2},\mu ^{n-2},..,\lambda ^{2},\mu ^{2},\lambda
^{1},\mu ^{1})d\mu ^{n-1}d\mu ^{n-2}..d\mu ^{2}d\mu ^{1}  \notag \\
&=&\dint \dint ..\dint \dint \mathcal{F}^{d}(\phi \ast _{c}F)(\xi
^{n},\lambda ^{n-1},\mu ^{n-1},\lambda ^{n-2},\mu ^{n-2},..,\lambda ^{2},\mu
^{2},\lambda ^{1},\mu ^{1})d\mu ^{n-1}d\mu ^{n-2}..d\mu ^{2}d\mu ^{1}  \notag
\\
&=&\dint \dint ..\dint \dint \mathcal{F}^{d}F(\xi ^{n},\lambda ^{n-1},\mu
^{n-1},\lambda ^{n-2},\mu ^{n-2},..,\lambda ^{2},\mu ^{2},\lambda ^{1},\mu
^{1})\mathcal{F}^{d}\phi (\xi ^{n},\lambda ^{n-1},\lambda ^{n-2},..,\lambda
^{2},\lambda ^{1})  \notag \\
&&d\mu ^{n-1}d\mu ^{n-2}..d\mu ^{2}d\mu ^{1}  \notag \\
&=&\mathcal{F}^{d}F(\xi ^{n},\lambda ^{n-1},0,\lambda ^{n-2},0,..,\lambda
^{2},0,\lambda ^{1},0)\mathcal{F}^{d}\phi (\xi ^{n},\lambda ^{n-1},\lambda
^{n-2},..,\lambda ^{2},\lambda ^{1})  \notag
\end{eqnarray}

\bigskip \textbf{Plancherels theorem 3.2}. For any function $f\in L^{1}(N),$
we have%
\begin{eqnarray}
&&\dint_{N}\left\vert \mathcal{F}^{d}f(\xi
^{n},X^{n-1},X^{n-2},..,X^{2},X^{1})\right\vert
^{2}dXdX^{n-1}dX^{n-2}..dX^{2}dX^{1}  \notag \\
&&\dint_{N}\left\vert \mathcal{F}^{d}f(\xi ^{n},\lambda ^{n-1},\lambda
^{n-2},..,\lambda ^{2},\lambda ^{1})\right\vert ^{2}d\xi ^{n}d\lambda
^{n-1}d\lambda ^{n-2}..d\lambda ^{2}d\lambda ^{1}
\end{eqnarray}

\bigskip \textit{Proof: Let }$f$ $\in L^{1}(K_{1})$ be the function defined
as 
\begin{eqnarray*}
&&\widetilde{\overset{\vee }{f}}%
(X^{n},U^{n-1},X^{n-1},U^{n-2},X^{n-2},..,U^{2},X^{2},U^{1},X^{1}) \\
&=&\overset{\vee }{f}[(\rho _{2}(X^{1})((\rho _{3}(X^{2})(..((\rho
_{n}(X^{n-1})(X^{n})),U^{n-1}+X^{n-1}),U^{n-2}+X^{n-2})..)),U^{2}+X^{2})),U^{1}+X^{1})]
\\
&=&f[(\rho _{2}(X^{1})((\rho _{3}(X^{2})(..((\rho
_{n}(X^{n-1})(X^{n})),U^{n-1}+X^{i-1}),U^{n-2}+X^{n-2})..)),U^{2}+X^{2})),U^{1}+X^{1})^{-1}]
\end{eqnarray*}%
then we have

\begin{eqnarray*}
&&\phi \ast \widetilde{\overset{\vee }{\phi }}(0,0,0,0,0,..0,0,0,0) \\
&=&\dint\limits_{N}\widetilde{\overset{\vee }{\phi }}\left[
(X^{n},X^{i-1},X^{i-2},..,X^{2},X^{1})^{-1}(0,0,0,0,0,..0,0,0,0)\right] \phi
(X^{n},X^{i-1},X^{i-2},..,X^{2},X^{1}) \\
&&dX^{n}dX^{n-1}dX^{n-2}..dX^{2}dX^{1} \\
&=&\dint\limits_{N}\widetilde{\overset{\vee }{\phi }}(\rho _{2}(X^{1})[(\rho
_{3}(X^{2})(..(\rho _{n-1}(X^{n-2})((\rho
_{n}(X^{n-1})(-X^{n})),0,-X^{i-1}),0,-X^{i-2},..,)0,-X^{2},) \\
&&0,-X^{1})]\phi
(X^{n},X^{i-1},X^{i-2},..,X^{2},X^{1})dX^{n}dX^{n-1}dX^{n-2}..dX^{2}dX^{1} \\
&=&\dint\limits_{N}\widetilde{\overset{\vee }{\phi }}\left[ (\rho
_{2}(X^{1})[(\rho _{3}(X^{2})(..(\rho _{n-1}(X^{n-2})((\rho
_{n}(X^{n-1})(-X^{n})),-X^{i-1}),-X^{i-2},..,),-X^{2},),-X^{1})\right] \\
&&\phi
(X^{n},X^{i-1},X^{i-2},..,X^{2},X^{1})dX^{n}dX^{n-1}dX^{n-2}..dX^{2}dX^{1} \\
&=&\dint\limits_{N}\overset{\vee }{\phi }\left[
(X^{n},X^{i-1},X^{i-2},..,X^{2},X^{1})^{-1}\right] \phi
(X^{n},X^{i-1},X^{i-2},..,X^{2},X^{1})dX^{n}dX^{n-1}dX^{n-2}..dX^{2}dX^{1} \\
&=&\dint\limits_{N}\overline{\phi (X^{n},X^{i-1},X^{i-2},..,X^{2},X^{1})}%
\phi
(X^{n},X^{i-1},X^{i-2},..,X^{2},X^{1})dX^{n}dX^{n-1}dX^{n-2}..dX^{2}dX^{1} \\
&=&\int\limits_{N}\left\vert \phi
(X^{n},X^{i-1},X^{i-2},..,X^{2},X^{1})\right\vert
^{2}dX^{n}dX^{n-1}dX^{n-2}..dX^{2}dX^{1}
\end{eqnarray*}

\bigskip In other hand we have%
\begin{eqnarray*}
&&\phi \ast \widetilde{\overset{\vee }{\phi }}(0,0,0,0,0,..0,0,0,0) \\
&=&\dint \mathcal{F}^{d}(\phi \ast \widetilde{\overset{\vee }{\phi }})(\xi
^{n},\lambda ^{n-1},\mu ^{n-1},\lambda ^{n-2},\mu ^{n-2},..,\lambda ^{2},\mu
^{2},\lambda ^{1},\mu ^{1}) \\
&&d\xi ^{n}d\lambda ^{n-1}d\mu ^{n-1}d\mu ^{n-2}d\lambda ^{n-2}..d\lambda
^{2}d\mu ^{2}d\lambda ^{1}d\mu ^{1} \\
&=&\dint \mathcal{F}^{d}(\phi \ast \widetilde{\overset{\vee }{\phi }})(\xi
^{n},\lambda ^{n-1},0,\lambda ^{n-2},0,..,\lambda ^{2},0,\lambda ^{1},0)d\xi
^{n}d\lambda ^{n-1}d\lambda ^{n-2}..d\lambda ^{2}d\lambda ^{1} \\
&=&\dint \mathcal{F}^{d}(\phi \ast _{c}\widetilde{\overset{\vee }{\phi }}%
)(\xi ^{n},\lambda ^{n-1},0,\lambda ^{n-2},0,..,\lambda ^{2},0,\lambda
^{1},0)d\xi ^{n}d\lambda ^{n-1}d\lambda ^{n-2}..d\lambda ^{2}d\lambda ^{1} \\
&=&\dint \mathcal{F}^{d}\widetilde{\overset{\vee }{\phi }}(\xi ^{n},\lambda
^{n-1},0,\lambda ^{n-2},0,..,\lambda ^{2},0,\lambda ^{1},0)\mathcal{F}%
^{d}\phi (\xi ^{n},\lambda ^{n-1},\lambda ^{n-2},..,\lambda ^{2},\lambda
^{1}) \\
&&d\xi ^{n}d\lambda ^{n-1}d\lambda ^{n-2}..d\lambda ^{2}d\lambda ^{1} \\
&=&\dint \overline{\mathcal{F}^{d}\phi (\xi ^{n},\lambda ^{n-1},\lambda
^{n-2},..,\lambda ^{2},\lambda ^{1})}\mathcal{F}^{d}\phi (\xi ^{n},\lambda
^{n-1},\lambda ^{n-2},..,\lambda ^{2},\lambda ^{1})d\xi ^{n}d\lambda
^{n-1}d\lambda ^{n-2}..d\lambda ^{2}d\lambda ^{1} \\
&=&\dint \left\vert \mathcal{F}^{d}\phi (\xi ^{n},\lambda ^{n-1},\lambda
^{n-2},..,\lambda ^{2},\lambda ^{1})\right\vert ^{2}d\xi ^{n}d\lambda
^{n-1}d\lambda ^{n-2}..d\lambda ^{2}d\lambda ^{1}
\end{eqnarray*}

Hence the proof of the theorem.

\subsection{Division of Distributions on $N$.}

\textbf{3.1.1.} Let $\mathcal{S}(N)$ be the Schwartz space of $N$ which%
\hspace{0.05in}can be considered as the Schwartz space of $\mathcal{S}(M),$
let $\mathcal{S}^{\prime }(N)$ be the space of all tempered distributions on 
$N$ and let $\mathcal{S}^{\prime }(M)$ be the space of all tempered
distributions on $M.$ If we consider the group $N$ as a subgroup of $K_{1}$,
then $\widetilde{f}\in \mathcal{S}(N)$ for $%
U=(U^{n-1},U^{n-2},..,U^{2},U^{1})$ is fixed, and if we consider $M$ as a
subgroup of $K_{1},$ then $\widetilde{f}\in \mathcal{S}(M)$ for $%
X=(X^{n-1},X^{n-2},..,X^{2},X^{1})$ is fixed. This being so; denote by $%
\mathcal{S}_{E}(K_{1})$ the space of all functions $\Phi
(X^{n},U^{n-1},X^{n-1},U^{n-2},X^{n-2},..,U^{2},X^{2},U^{1},X^{1})\in
C^{\infty }(K_{1})$ such that $\Phi
(X^{n},U^{n-1},X^{n-1},U^{n-2},X^{n-2},..,U^{2},X^{2},U^{1},X^{1})\in 
\mathcal{S}(N)$ for $X=(X^{n-1},X^{n-2},..,X^{2},X^{1})$ is fixed, and $\Phi
(X^{n},U^{n-1},X^{n-1},U^{n-2},X^{n-2},..,U^{2},X^{2},U^{1},X^{1})\in 
\mathcal{S}(M)$ for $U=(U^{n-1},U^{n-2},..,U^{2},U^{1})$ is fixed. We equip $%
\mathcal{S}_{E}(K_{1})$ with the natural topology defined by the seminormas: 
\begin{equation*}
\Phi \rightarrow \underset{(X^{n},U)\in M}{\text{sup}}\left\vert
Q(Y)P(D)\Phi
(X^{n},U^{n-1},X^{n-1},U^{n-2},X^{n-2},..,U^{2},X^{2},U^{1},X^{1})\right%
\vert \qquad X\text{ }fixed
\end{equation*}%
\begin{equation*}
\Phi \rightarrow \underset{(X^{n},X)\in N\text{ }}{\text{sup}}\left\vert
R(Y)S(D)\Phi
(X^{n},U^{n-1},X^{n-1},U^{n-2},X^{n-2},..,U^{2},X^{2},U^{1},X^{1})\right%
\vert \text{ \qquad }U\text{ }fixed
\end{equation*}%
where $%
(X^{n},X)=(X^{n},X^{n-1},X^{n-2},..,X^{2},X^{1}),(X^{n},U)=(X^{n},U^{n-1},U^{n-2},..,U^{2},U^{1}), 
$ and $P,$ $Q,$ $R$ and $S$ run over the family of all complex polynomials
in $2d-n$ variables. Let $\mathcal{S}_{E}^{I}(K_{1})$ be the subspace of all
functions $\psi \in \mathcal{S}_{E}(K_{1}),$ which are invariant in sense $%
(13)$, then we have the following result.\bigskip \bigskip

\textbf{Lemma 3.1.1.} \textit{Let }$u\in \mathcal{U}$\textit{\ and }$Q_{u}$%
\textit{\ be the invariant differential operator on the group }$M,$ which is
associated to $u,$ \textit{acts on the variables }$%
(X^{n},U^{n-1},U^{n-2},..,U^{2},U^{1})\in M,$ then we have

$(i)$ \textit{The mapping }$f\mapsto \widetilde{f}$\textit{\ is a
topological isomorphism of }$\mathcal{S}(N)$\textit{\ onto }$\mathcal{S}%
_{E}^{I}(K_{1})$\textit{.}

$(ii)$ \textit{The mapping }$\Phi \mapsto Q_{u}\Phi $\textit{\ is a
topological isomorphism of }$\mathcal{S}_{E}^{I}(K_{1})$\textit{\ onto its
image}

\textit{Proof}\textbf{: }$(i)$ In fact $\sim $ is continuous and the
restriction mapping $\Phi \mapsto R\Phi $ on $N$ is continuous from $%
\mathcal{S}_{E}^{I}(K_{1})$ into $\mathcal{S}(N)$ that satisfies $R\circ
\sim =Id_{\mathcal{S}(N)}$ and $\sim \circ R=Id_{\mathcal{S}%
_{E}^{I}(K_{1})}, $ where $Id_{\mathcal{S}(N)}$ $(resp$. $Id_{\mathcal{S}%
_{E}^{I}(K_{1})})$ is the identity mapping of $\mathcal{S}(N)$ $(resp$. $%
\mathcal{S}_{E}^{I}(K_{1}))$ and $N$ is considered as a subgroup of $K_{1}.$

To prove$(ii)$ we refer to$[16,$ $P.313-315]$ and his famous result that is:
"\textit{Any invariant differential operator on }$M,$\textit{\ is a
topological isomorphism of \ }$\mathcal{S}(M)$\textit{\ onto its image}"
From this result, we obtain that 
\begin{equation}
D_{u}:\mathcal{S}_{E}(K_{1})\rightarrow \mathcal{S}_{E}(K_{1})
\end{equation}%
is a topological isomorphism and its restriction on $\mathcal{S}%
_{E}^{I}(K_{1})$ is a topological isomorphism of $\mathcal{S}_{E}^{I}(K_{1})$
onto its image. Hence the theorem is proved. \ 

In the following we will prove that every invariant differential operator on 
$N$ has a tempered fundamental solution. As in the introduction, we will
consider the two invariant differential operators $P_{u}$ and $Q_{u}$, the
first on the group $N$ $,$ and the second on the abelian vector group $M.$
Our main result is:

\textbf{Theorem 3.1.1}. \textit{Every nonzero invariant differential
operator on }$N$ \textit{has a tempered fundamental solution}

\textit{Proof}\textbf{:}\textit{\ For every function} $\psi \in C^{\infty
}(K_{1})$ \textit{invariant in sense} $(11)$ \textit{and for every }$u\in $ 
\textit{$\mathcal{U}$}, \textit{we have} 
\begin{eqnarray*}
&&P_{u}\psi
(X^{n},U^{n-1},X^{n-1},U^{n-2},X^{n-2},..,U^{2},X^{2},U^{1},X^{1}) \\
&=&u\text{ }\ast \psi
(X^{n},U^{n-1},X^{n-1},U^{n-2},X^{n-2},..,U^{2},X^{2},U^{1},X^{1}) \\
&=&u\text{ }\ast _{c}\psi
(X^{n},U^{n-1},X^{n-1},U^{n-2},X^{n-2},..,U^{2},X^{2},U^{1},X^{1}) \\
&=&Q_{u}\psi
(X^{n},U^{n-1},X^{n-1},U^{n-2},X^{n-2},..,U^{2},X^{2},U^{1},X^{1})
\end{eqnarray*}%
\textit{for every} $%
(X^{n},U^{n-1},X^{n-1},U^{n-2},X^{n-2},..,U^{2},X^{2},U^{1},X^{1})$ $\in
K_{1}$, \textit{where} $\ast $\textit{\ signifies the convolution product on}
$N$ \textit{with respect the variables} ($%
X^{n},X^{n-1},X^{n-2},..,X^{2},X^{1})$\ \textit{and} $\ast _{c}$\textit{%
signifies the commutative convolution product on} $M$ \textit{with respect
the variables (}$X^{n},U^{n-1},U^{n-2},..,U^{2},U^{1}).$

\textit{Proof}\textbf{:}\textit{\ }\textbf{\ }In fact we have.\ By Lemma $%
3.1.1,$ the mapping $\psi \mapsto Q_{u}\psi $ is a topological isomorphism
of $\mathcal{S}_{E}^{I}(K_{1})$ onto its image, then the mapping $\psi
\mapsto P_{u}\psi $ \ is a topological isomorphism of $\mathcal{S}%
_{E}^{I}(K_{1})$ onto its image.

Since 
\begin{eqnarray}
&&R(P_{u}\psi
)(X^{n},U^{n-1},X^{n-1},U^{n-2},X^{n-2},..,U^{2},X^{2},U^{1},X^{1})  \notag
\\
&=&P_{u}(R\psi
)(X^{n},U^{n-1},X^{n-1},U^{n-2},X^{n-2},..,U^{2},X^{2},U^{1},X^{1})
\end{eqnarray}%
\ so the following diagram is commutative:\bigskip

$\bigskip \qquad \qquad \qquad \qquad \qquad \mathcal{S}_{E}^{I}(K_{1})$ $%
\qquad $\qquad $\underrightarrow{\underset{}{P_{u}}}\qquad $\qquad $P_{u}%
\mathcal{S}_{E}^{I}(K_{1})\bigskip $

$\qquad \qquad \qquad \qquad $\qquad $\quad \sim \uparrow \downarrow R$
\qquad \qquad \qquad $\qquad $\qquad $\downarrow R\bigskip $

$\qquad \qquad \qquad $\qquad $\qquad \mathcal{S}(N)$ \qquad $\qquad 
\underrightarrow{\underset{}{Q_{u}}}$\qquad $\qquad P_{u}\mathcal{S}%
(N)\bigskip $

Hence the mapping $\psi \mapsto P_{u}\psi $ \ is a topological isomorphism
of $\mathcal{S}(N)$ onto its image. So the transpose $^{t}P_{u}$ of $P_{u}$
is a continuous mapping of $\mathcal{S}^{\prime }(N)$ onto $\mathcal{S}%
^{\prime }(N).$ This means that for every tempered distribution $T$ on $N$
there is a tempered distribution $E$ on $N$ such that 
\begin{equation}
P_{u}F=T
\end{equation}

Indeed the Dirac measure $\delta $ belongs to $\mathcal{S}^{\prime }(N).$

\bigskip If we use Atiyah method $[1],$ and $[9],$ we can obtain the
following result "\textit{Every invariant differential operator on }$N$ 
\textit{which is not identically }$0$\textit{\ has a tempered fundamental
solution".}

\section{\protect\bigskip Fourier Analysis on $AN.$}

\bigskip \textbf{4.1.} Let $G=SL(n,\mathbb{R})$\ be the real semi-simple Lie
group and let $G=KAN$ \ be the Iwasawa decomposition of $G$, where $K=SO(n,%
\mathbb{R}),$and%
\begin{equation}
N=\left( 
\begin{array}{ccccc}
1 & \ast & . & . & \ast \\ 
0 & 1 & \ast & . & \ast \\ 
. & . & . & . & \ast \\ 
. & . & . & . & . \\ 
0 & 0 & . & 0 & 1%
\end{array}%
\right) ,
\end{equation}

\begin{equation}
A=\left( 
\begin{array}{ccccc}
a_{1} & 0 & 0 & . & 0 \\ 
0 & a_{2} & 0 & . & 0 \\ 
. & . & . & . & . \\ 
. & . & . & . & . \\ 
0 & 0 & . & 0 & a_{n}%
\end{array}%
\right)
\end{equation}%
where $a_{1}.a_{2}....a_{n}=1$ and $a_{i}\in \mathbb{R}_{+}^{\ast }.$ The
product $AN$ is a closed subgroup of $G$ and is isomorphic (algebraically
and topologically) to the semi-direct product of $A$ and $N$ with $N$ normal
in $AN.$

Then the group $AN$ is nothing but the group $S=$ $N\rtimes $ $_{\rho }A.$
So the product of two elements $X$ and$Y$ by%
\begin{equation}
(x,\text{ }a)(m,\text{ }b)=(x.\rho (a)y,\text{ }a.b)
\end{equation}%
for\ any $X=(x,a_{1},a_{2},..,a_{n})\in S$ and $Y$ $%
=(m,b_{1},b_{2},..,b_{n})\in S.$ Let $dnda=dmda_{1}da_{2}..da_{n-1}$ be the
right haar measure on $S$ and let $L^{2}(S)$ be the Hilbert space of the
group $S.$ Let $L^{1}(S)$ be the Banach algebra that consists of all complex
valued functions on the group $S$, which are integrable with respect to the
Haar measure of $S$ and multiplication is defined by convolution on $S$ as

\begin{equation}
g\ast f=\int\limits_{S}f((m,b)^{-1}(n,a))g(m,b)dmdb
\end{equation}%
where $dmdb=dmdb_{1}db_{2}..db_{n-1}$ is the right Haar measure on $S=$ $%
N\rtimes $ $_{\rho }A.$

\bigskip In the following we prove the Plancherel theorem. Therefore let $%
T=N\times A$ be the Lie group of direct product of the two Lie groups $N$
and $A,$ and let $H=N\times A\times A$ the Lie group, with law%
\begin{equation}
(n,t,r)(m,s,q)=(n\rho (r)m,ts,rq)
\end{equation}%
for all $(n,t,r)\in H$ and $(m,s,q)\in H.$ In this case the group $S$ can be
identified with the closed subgroup $N\times \left\{ 0\right\} \times _{\rho
}A$ of $H$ and $T$ with the subgroup $N\times A\times \left\{ 0\right\} $of $%
H$

\ \textbf{Definition 4.1.}\textit{\ For every functon }$f$ defined on $S$, 
\textit{one can define a function} $\widetilde{f}$ on $L$ \textit{as follows:%
} 
\begin{equation}
\widetilde{f}(n,a,b)=f(\rho (a)n,ab)
\end{equation}%
\textit{for all} $(n,a,b)\in H.$ \textit{So every function} $\psi (n,a)$ 
\textit{on} $S$\textit{\ extends uniquely as an invariant function} $%
\widetilde{\psi }(n,$ $b,$ $a)$ \textit{on} $L.$

\ \textbf{Remark 4.1. }\textit{The function} $\widetilde{f}$ \textit{is
invariant in the following sense:} 
\begin{equation}
\widetilde{f}(\rho (s)n,as^{-1},bs)=\widetilde{f}(n,a,b)
\end{equation}%
\textit{for any} $(n,a,b)\in H$ \textit{and} $s\in H.$

\textbf{Lemma 4.1.}\textit{\ For every function} $f\in L^{1}(S)$ \textit{and
for every }$g\mathcal{\in }$ $L^{1}(S)$, \textit{we have} 
\begin{equation}
g\ast \widetilde{f}(n,a,b)=g\ast _{c}\widetilde{f}(n,a,b)
\end{equation}%
\begin{equation}
\int\limits_{\mathbb{R}^{n-1}}\mathcal{F}^{d}\mathcal{(}g\ast \widetilde{f}%
)(\lambda ,\mu ,\nu )d\nu =\mathcal{F}^{d}\widetilde{f}(\lambda ,\mu ,0)%
\mathcal{F}_{A}^{n-1}g(\lambda ,\mu )
\end{equation}%
\textit{for every} $(n,a,b)$ $\in H$, \textit{where} $\ast $\textit{\
signifies the convolution product on} $S$\ \textit{with respect the variables%
} $(n,b)$\ \textit{and} $\ast _{c}$\textit{signifies the commutative
convolution product on} $B$ \textit{with respect the variables }$(n,a),$ 
\textit{where }$\mathcal{F}_{A}^{n-1}$ \textit{is the Fourier transform on }$%
A$ \textit{and }$\mathcal{F}^{d}$ \textit{is the Fourier transform on} $N$

\bigskip \textit{Proof: }\textbf{\ }In fact we have\ 
\begin{eqnarray}
&&g\ast \widetilde{f}(n,a,b)=\int\limits_{G}\widetilde{f}%
((m,c)^{-1}(n,a,b))g(m,s)dmdc  \notag \\
&=&\int\limits_{S}\widetilde{f}\left[ (\rho (c^{-1})(m^{-1}),c^{-1})(n,a,b)%
\right] u(m,s)dmdc  \notag \\
&=&\int\limits_{S}\widetilde{f}\left[ \rho (c^{-1})(m^{-1}n),a,bc^{-1}\right]
g(m,c)dmdc  \notag \\
&=&\int\limits_{S}\widetilde{f}\left[ m^{-1}n,ac^{-1},b\right]
g(m,c)dmdc=g\ast _{c}\widetilde{f}(n,a,b)
\end{eqnarray}%
and so%
\begin{eqnarray}
&&\int\limits_{\mathbb{R}^{n-1}}\mathcal{F}^{d}\mathcal{F}_{A}^{n-1}\mathcal{%
(}g\ast \widetilde{f})(\lambda ,\mu ,\nu )d\nu  \notag \\
&=&\int\limits_{\mathbb{R}^{n-1}}\mathcal{F}^{d}\mathcal{F}_{A}^{n-1}%
\mathcal{(}g\ast _{c}\widetilde{f})(\lambda ,\mu ,\nu )d\nu  \notag \\
&=&\mathcal{F}^{d}\mathcal{F}_{A}^{n-1}\widetilde{f}(\lambda ,\mu ,0)%
\mathcal{F}_{A}^{n-1}g(\lambda ,\mu )
\end{eqnarray}

\textbf{Plancherel} \textbf{Theorem 4.1. }\textit{For any function }$\Psi
\in L^{1}(S),$ \textit{we have }%
\begin{equation}
\int\limits_{\mathbb{R}^{d}}\int\limits_{\mathbb{R}^{n-1}}\left\vert 
\mathcal{F}^{d}\mathcal{F}_{A}^{n-1}\Psi (\lambda ,\mu )\right\vert
^{2}d\lambda =\int\limits_{AN}\left\vert \Psi (X,a)\right\vert ^{2}\text{ }dX
\end{equation}

\textit{Proof}: Let $\widetilde{\overset{\vee }{\Psi }}$ the function
defined as 
\begin{equation}
\widetilde{\overset{\vee }{\Psi }}(n,a,b)=\overset{\vee }{\Psi }(an,ab)=%
\overline{\Psi \lbrack (an,ab)^{-1}]}
\end{equation}

Then we obtain

\begin{eqnarray}
&&\Psi \mathcal{\ast }\widetilde{\overset{\vee }{\Psi }}(0,0,0)  \notag \\
&=&\int\limits_{NA}\widetilde{\overset{\vee }{\Psi }}((n,a)^{-1}(0,0,0))\Psi
((n,a)dnda  \notag \\
&=&\int\limits_{NA}\widetilde{\overset{\vee }{\Psi }}%
(a^{-1}.n,a^{-1})(0,0,0))\Psi ((n,a)dnda=\int\limits_{NA}\widetilde{\overset{%
\vee }{\Psi }}(a^{-1}.n,0,a^{-1})\Psi ((n,a)dnda  \notag \\
&=&\int\limits_{NA}\widetilde{\overset{\vee }{\Psi }}((-X^{1})(X),0,-X^{1})%
\Psi ((n,a)dnda=\int\limits_{NA}\overset{\vee }{\Psi }(a^{-1}.n,a^{-1})\Psi
((n,a)dnda  \notag \\
&=&\int\limits_{NA}\overline{\overset{\vee }{\Psi }(a^{-1}.n,a)^{-1}}\Psi
((n,a)dnda=\int\limits_{\mathbb{R}^{d}}\left\vert \Psi (n,a)\right\vert
^{2}dnda
\end{eqnarray}%
and by lemma 4.1, we get%
\begin{eqnarray*}
&&\Psi \mathcal{\ast }\widetilde{\overset{\vee }{\Psi }}(0,0,0) \\
&=&\int\limits_{N}\int\limits_{A}\mathcal{F}^{d}\mathcal{F}_{A}^{n-1}(\Psi 
\mathcal{\ast }\widetilde{\overset{\vee }{\Psi }})(\xi ,\lambda ,\mu )d\xi
d\lambda d\mu \\
&=&\int\limits_{N}\int\limits_{A}\mathcal{F}^{d}\mathcal{F}_{A}^{n-1}(\Psi 
\mathcal{\ast }_{c}\widetilde{\overset{\vee }{\Psi }})(\xi ,\lambda ,\mu
)d\xi d\lambda d\mu \\
&=&\int\limits_{N}\int\limits_{A}\mathcal{F}^{d}\mathcal{F}_{A}^{n-1}\text{ }%
\widetilde{\overset{\vee }{\Psi }}(\xi ,\lambda ,\mu )\mathcal{F}^{d}%
\mathcal{F}_{A}^{n-1}\Psi (\xi ,\lambda )d\xi d\lambda d\mu \\
&=&\int\limits_{N}\int\limits_{A}\int\limits_{N}\int\limits_{A}\int%
\limits_{N}\int\limits_{A}\widetilde{\overset{\vee }{\Psi }}%
(n,a,0)e^{-i\langle \xi ,n\rangle }e^{-i\langle \lambda ,a\rangle }\Psi
(m,b)e^{-i\langle \xi ,m\rangle }e^{-i\langle \lambda ,b\rangle } \\
&&dndadmdbd\xi d\lambda \\
&=&\int\limits_{N}\int\limits_{A}\int\limits_{N}\int\limits_{A}\int%
\limits_{N}\int\limits_{A}\overset{\vee }{\Psi }(an,a)e^{-i\langle \xi
,n\rangle }e^{-i\langle \lambda ,a\rangle }\Psi (m,b)e^{-i\langle \xi
,m\rangle }e^{-i\langle \lambda ,b\rangle } \\
&&dndadmdbd\xi d\lambda \\
&=&\int\limits_{N}\int\limits_{A}\int\limits_{N}\int\limits_{A}\int%
\limits_{N}\int\limits_{A}\overline{\Psi (an,a)^{-1}}e^{-i\langle \xi
,n\rangle }e^{-i\langle \lambda ,a\rangle }\Psi (m,b)e^{-i\langle \xi
,m\rangle }e^{-i\langle \lambda ,b\rangle } \\
&&dndadmdbd\xi d\lambda \\
&=&\int\limits_{N}\int\limits_{A}\int\limits_{N}\int\limits_{A}\int%
\limits_{N}\int\limits_{A}\overline{\Psi (n^{-1},-a)e^{-i\langle \xi
,n\rangle }e^{-i\langle \lambda ,a\rangle }}\Psi (m,b)e^{-i\langle \xi
,m\rangle }e^{-i\langle \lambda ,b\rangle } \\
&&dndadmdbd\xi d\lambda \\
&=&\int\limits_{N}\int\limits_{A}\int\limits_{N}\int\limits_{A}\int%
\limits_{N}\int\limits_{A}\overline{\Psi (n,a)e^{-i\langle \xi ,n\rangle
}e^{-i\langle \lambda ,a\rangle }}\Psi (m,b)e^{-i\langle \xi ,m\rangle
}e^{-i\langle \lambda ,b\rangle } \\
&=&\int\limits_{N}\int\limits_{A}\overline{\mathcal{F}^{d}\mathcal{F}%
_{A}^{n-1}(\Psi )(\xi ,\lambda )}\mathcal{F}^{d}\mathcal{F}_{A}^{n-1}(\Psi
)(\xi ,\lambda )d\xi d\lambda \\
&=&\int\limits_{N}\int\limits_{A}\left\vert \mathcal{F}^{d}\mathcal{F}%
_{A}^{n-1}(\Psi )(\xi ,\lambda \right\vert ^{2}d\xi d\lambda
\end{eqnarray*}%
where $\langle \xi ,n\rangle =\dsum\limits_{i=1}^{d}\xi _{i}n_{i},$ and $%
\langle \lambda ,a\rangle =\dsum\limits_{j=1}^{n-1}\lambda _{j}n_{j}.$

\textbf{4.2.} Let $\mathcal{S}(N\rtimes A)$ be the Schwartz space of $%
S=N\rtimes A$ which\hspace{0.05in}can be considered as the Schwartz space of 
$\mathcal{S}(T=N\times A),$ let $\mathcal{S}^{\prime }(N\rtimes A)$ be the
space of all tempered distributions on $N\rtimes A$ and let $\mathcal{S}%
^{\prime }(T)$ be the space of all tempered distributions on $C.$ If we
consider the group $S$ as a subgroup of $H$, then $\widetilde{f}\in \mathcal{%
S}(S)$ for $a$ is fixed, and if we consider $T$ as a subgroup of $H,$ then $%
\widetilde{f}\in \mathcal{S}(T)$ for $b$ is fixed. This being so; denote by $%
\mathcal{S}_{E}(H)$ the space of all functions $\Phi (x,a,b)\in C^{\infty
}(H)$ such that $\Phi (x,a,b))\in \mathcal{S}(N)$ for $a$ is fixed, and $%
\Phi (x,a,b)\in \mathcal{S}(T)$ for $b$ is fixed. We equip$\mathcal{S}%
_{E}(H) $ with the natural topology defined by the seminormas: 
\begin{equation*}
\Phi \rightarrow \underset{(x,a)\in T}{\text{sup}}\left\vert Q(Y)P(D)\Phi
(x,a,b)\right\vert \qquad b\text{ }fixed
\end{equation*}%
\begin{equation*}
\Phi \rightarrow \underset{(x,b)\in N\text{ }}{\text{sup}}\left\vert
R(Y)S(D)\Phi (x,a,b)\right\vert \text{ \qquad }a\text{ }fixed
\end{equation*}%
where $P,$ $Q,$ $R$ and $S$ run over the family of all complex polynomials
in $2d-1$ variables. Let $\mathcal{S}_{E}^{I}(H)$ be the subspace of all
functions $\psi \in \mathcal{S}_{E}(H),$ which are invariant in sense $(13)$%
, then we have the following result.\bigskip \bigskip

\textbf{Proposition 4.1.} \textit{Let }$u\in \mathcal{U}$\textit{\ and }$%
D_{u}$\textit{\ be the invariant differential operator on the group }$T,$
which is associated to $u,$ \textit{acts on the variables }$(x,a)\in T,$
then we have

$(i)$ \textit{The mapping }$f\mapsto \widetilde{f}$\textit{\ is a
topological isomorphism of }$\mathcal{S}(N)$\textit{\ onto }$\mathcal{S}%
_{E}^{I}(H)$\textit{.}

$(ii)$ \textit{The mapping }$\Phi \mapsto Q_{u}\Phi $\textit{\ is a
topological isomorphism of }$\mathcal{S}_{E}^{I}(H)$\textit{\ onto its image}

\textit{Proof}\textbf{: }$(i)$ In fact $\sim $ is continuous and the
restriction mapping $\Phi \mapsto R\Phi $ on $S$ is continuous from $%
\mathcal{S}_{E}^{I}(H)$ into $\mathcal{S}(S)$ that satisfies $R\circ \sim
=Id_{\mathcal{S}(N\rtimes A)}$ and $\sim \circ R=Id_{\mathcal{S}%
_{E}^{I}(H)}, $ where $Id_{\mathcal{S}(S)}$ $(resp$. $Id_{\mathcal{S}%
_{E}^{I}(H)})$ is the identity mapping of $\mathcal{S}(S)$ $(resp$. $%
\mathcal{S}_{E}^{I}(H))$ and $S$ is considered as a subgroup of $H.$

To prove$(ii)$ we refer to$[16,$ $P.313-315]$ and his famous result that is:
"\textit{Any invariant differential operator on }$T,$\textit{\ is a
topological isomorphism of \ }$\mathcal{S}(T)$\textit{\ onto its image}"
From this result, we obtain that 
\begin{equation}
D_{u}:\mathcal{S}_{E}(H)\rightarrow \mathcal{S}_{E}(H)
\end{equation}%
is a topological isomorphism and its restriction on $\mathcal{S}_{E}^{I}(H)$
is a topological isomorphism of $\mathcal{S}_{E}^{I}(H)$ onto its image.
Hence the proposition is proved. \ 

In the following we will prove that every invariant differential operator on 
$S=N\rtimes A$ has a tempered fundamental solution. We consider the two
invariant differential operators $P_{u}$ and $Q_{u}$, the first on the group 
$S$ $,$ and the second on the abelian vector group $T.$ Our main result is:

\textbf{Existence theorem 4.2}. \textit{Every nonzero invariant differential
operator on }$N\rtimes A$ \textit{has a tempered fundamental solution}

\textit{Proof}\textbf{:}\textit{\ For every function} $\psi \in C^{\infty
}(H)$ \textit{invariant in sense} $(11)$ \textit{and for every }$u\in $ 
\textit{$\mathcal{U}$}, \textit{we have} 
\begin{eqnarray*}
&&P_{u}\psi (x,a,b)=u\text{ }\ast \psi (x,a,b) \\
&=&u\text{ }\ast _{c}\psi (x,a,b)=Q_{u}\psi (x,a,b)
\end{eqnarray*}%
\textit{for every} $(x,a,b)$ $\in H$, \textit{where} $\ast $\textit{\
signifies the convolution product on} $S$ \textit{with respect the variables}
$(x,b)$\ \textit{and} $\ast _{c}$\textit{signifies the commutative
convolution product on} $T$ \textit{with respect the variables }$(x,a).$

\textit{Proof}\textbf{:}\textit{\ }\textbf{\ }In fact we have by proposition 
$4.1,$ the mapping $\psi \mapsto Q_{u}\psi $ is a topological isomorphism of 
$\mathcal{S}_{E}^{I}(H)$ onto its image, then the mapping $\psi \mapsto
P_{u}\psi $ \ is a topological isomorphism of $\mathcal{S}_{E}^{I}(H)$ onto
its image.

Since 
\begin{equation}
R(P_{u}\psi )(x,a,b)=P_{u}(R\psi )(x,a,b)
\end{equation}%
\ so the following diagram is commutative:\bigskip

$\bigskip \qquad \qquad \qquad \qquad \qquad \mathcal{S}_{E}^{I}(H)$ $\qquad 
$\qquad $\underrightarrow{\underset{}{P_{u}}}\qquad $\qquad $P_{u}\mathcal{S}%
_{E}^{I}(H)\bigskip $

$\qquad \qquad \qquad \qquad $\qquad $\quad \sim \uparrow \downarrow R$
\qquad \qquad \qquad $\qquad $\qquad $\downarrow R\bigskip $

$\qquad \qquad \qquad $\qquad $\qquad \mathcal{S}(N\rtimes A)$ \qquad $%
\qquad \underrightarrow{\underset{}{Q_{u}}}$\qquad $\qquad P_{u}\mathcal{S}%
(N\rtimes A)\bigskip $

Hence the mapping $\psi \mapsto P_{u}\psi $ \ is a topological isomorphism
of $\mathcal{S}(N\rtimes A)$ onto its image. So the transpose $^{t}P_{u}$ of 
$P_{u}$ is a continuous mapping of $\mathcal{S}^{\prime }(N\rtimes A)$ onto $%
\mathcal{S}^{\prime }(N\rtimes A).$ This means that for every tempered
distribution $F$ on $N\rtimes A$ there is a tempered distribution $E$ on $%
N\rtimes A$ such that 
\begin{equation}
P_{u}E=F
\end{equation}%
Indeed the Dirac measure $\delta $ belongs to $\mathcal{S}^{\prime
}(N\rtimes A).$

\section{Ideals Algebra $L^{1}(N)$.}

\textbf{Lemma 5.1. }$(i)$\textit{The mapping }$\Gamma $ \textit{from }$%
\widetilde{L^{1}(N)}|_{M}$ \textit{to\ }$\widetilde{L^{1}(N)}|_{N\text{ }}$ 
\textit{defined by } 
\begin{eqnarray*}
&&\widetilde{\phi }|_{M}\ (X^{n},X^{n-1},0,X^{n-2},0..,X^{2},0,X^{1},0) \\
&\rightarrow &\Gamma (\widetilde{\phi }|_{M}\
)(X^{n},0,X^{n-1},0,X^{n-2},..,0,X^{2},0,X^{1}) \\
&=&\widetilde{\phi }|_{N}(X^{n},0,X^{n-1},0,X^{n-2},..,0,X^{2},0,X^{1})\ 
\end{eqnarray*}%
\textit{\ \ is a topological isomorphism}

\bigskip $(ii)$ \textit{For every }$\psi \in L^{1}(N)$ \textit{and} $\phi
\in L^{1}(N),$ \textit{we obtain} 
\begin{eqnarray*}
&&\Gamma (\psi \ast _{c}\widetilde{\phi }%
|_{M})(X^{n},0,X^{n-1},0,X^{n-2},..,0,X^{2},0,X^{1}) \\
&=&\psi \ast \widetilde{\phi }%
|_{N}(X^{n},0,X^{n-1},0,X^{n-2},..,0,X^{2},0,X^{1}) \\
&=&\psi \ast \widetilde{\phi }(X^{n},X^{n-1},X^{n-2},..,X^{2},X^{1})
\end{eqnarray*}%
\begin{eqnarray*}
&&(\psi \ast _{c}\widetilde{\phi }%
|_{M})(X^{n},X^{n-1},0,X^{n-2},0,..,X^{2},0,X^{1},0)\text{ \ \ \ \ \ \ \ \ \
\ \ \ \ \ } \\
&=&=\int\limits_{M}\widetilde{\phi }\left[
(X^{n}-Y^{n},X^{n-1}-Y^{n-1},0,X^{n-2}-Y^{n-2},0,..,X^{2}-Y^{2},0,X^{1}-Y^{1},0%
\right] \\
&&\psi
(Y^{n},Y^{n-1},Y^{n-2},..,Y^{2},Y^{1})dY^{n}dY^{n-1}dY^{n-2},..,dY^{2}dY^{1}%
\text{ \ \ \ \ \ \ \ \ \ \ \ \ \ \ \ \ \ \ \ \ \ \ \ \ \ \ \ \ \ \ \ \ \ \ \
\ \ \ }
\end{eqnarray*}

\textit{Proof}\textbf{: }$(i)$ The mapping $\Gamma $ is continuous and has
an inverse $\Gamma ^{-1}$ given by 
\begin{eqnarray*}
&&\widetilde{\phi }|_{N}\ (X^{n},0,X^{n-1},0,X^{n-2},..,0,X^{2},0,X^{1}) \\
&\rightarrow &\Gamma ^{-1}(\widetilde{\phi }|_{N}\
)(X^{n},X^{n-1},0,X^{n-2},0..,X^{2},0,X^{1},0) \\
&=&\widetilde{\phi }|_{M}(X^{n},X^{n-1},0,X^{n-2},0..,X^{2},0,X^{1},0)\ 
\end{eqnarray*}

$(ii)$ It is enough to see for every $\phi \in L^{1}(N)$ 
\begin{eqnarray*}
&&\Gamma (\psi \ast _{c}\widetilde{\phi }%
|_{M})(X^{n},0,X^{n-1},0,X^{n-2},..,0,X^{2},0,X^{1}) \\
&=&\int\limits_{M}\widetilde{\phi }|_{M}\left[
X^{n}-Y^{n},-Y^{n-1},X^{n-1},-Y^{n-2},X^{n-2},..,-Y^{2},X^{2},-Y^{1},X^{1}%
\right] \\
&&\psi
(Y^{n},Y^{n-1},Y^{n-2},..,Y^{2},Y^{1})dY^{n}dY^{n-1}dY^{n-2},..,dY^{2}dY^{1}
\\
&=&\psi \ast \phi (X^{n},X^{n-1},X^{n-2},..,X^{2},X^{1}),\text{ \ \ \ \ \ \
\ \ \ \ \ \ \ \ \ \ \ \ \ \ \ \ \ \ \ \ \ \ }\phi \in L^{1}(N)
\end{eqnarray*}

\bigskip If $I$ is a subalgebra of $L^{1}(N),$ we denote by $\widetilde{I}$
its image by the mapping $\thicksim $. Let $J=$ $\widetilde{I}\ |_{M}.$ Our
main result is:

\textbf{Theorem 5.1.} \textit{Let }$I$\textit{\ be a subalgebra of }$%
L^{1}(N),$\textit{\ then the following conditions are equivalents.}

$(i)$\textit{\ }$J=\widetilde{I}\ |_{M}$\textit{\ is an ideal in the Banach
algebra }$L^{1}(M).$

$(ii)$ $I$ \textit{is a left ideal in the Banach algebra }$L^{1}(N).$

\textit{Proof:}\textbf{\ }$(i)$ implies $(ii)$\ Let $I$ be a subspace of the
space $L^{1}(N)$ such that $J=\widetilde{I}|_{M}$ \ is an ideal in $%
L^{1}(M), $ then we have: 
\begin{eqnarray}
&&\psi \ast _{c}\widetilde{I}\
|_{M}(X^{n},X^{n-1},0,X^{n-2},0..,X^{2},0,X^{1},0)  \notag \\
&\subseteq &\widetilde{I}\ |_{M}(X^{n},X^{n-1},0,X^{n-2},0..,X^{2},0,X^{1},0)
\end{eqnarray}%
for any $\psi \in L^{1}(M)$ and $(x,x_{3},x_{2},x_{1})\in M$, where 
\begin{eqnarray*}
&&\psi \ast _{c}\widetilde{I}\
|_{M}(X^{n},X^{n-1},0,X^{n-2},0..,X^{2},0,X^{1},0) \\
&=&\left\{ 
\begin{array}{c}
\int\limits_{M}\widetilde{\phi }\left[
X^{n}-Y^{n},-Y^{n-1},X^{n-1},-Y^{n-2},X^{n-2},..,-Y^{2},X^{2},-Y^{1},X^{1}%
\right] \\ 
\psi (Y^{n},Y^{n-1},Y^{n-2},..,Y^{2},Y^{1} \\ 
dY^{n}dY^{n-1}dY^{n-2},..,dY^{2}dY^{1},\text{ \ \ \ \ }\phi \in L^{1}(N)%
\end{array}%
\right\}
\end{eqnarray*}

It shows that 
\begin{eqnarray}
&&\psi \ast _{c}\widetilde{\phi }\
|_{M}(X^{n},X^{n-1},0,X^{n-2},0..,X^{2},0,X^{1},0)  \notag \\
&\in &\widetilde{I}\ |_{M}(X^{n},X^{n-1},0,X^{n-2},0..,X^{2},0,X^{1},0)
\end{eqnarray}%
for any $\widetilde{\phi }\in \widetilde{I}.$ Then we get 
\begin{eqnarray}
&&\Gamma (\psi \ast _{c}\widetilde{\phi }%
|_{M})(X^{n},0,X^{n-1},0,X^{n-2},..,0,X^{2},0,X^{1})  \notag \\
&=&\psi \ast \widetilde{\phi }(X^{n},0,X^{n-1},0,X^{n-2},..,0,X^{2},0,X^{1})
\notag \\
&\in &\Gamma (\widetilde{I}\
|_{M})(X^{n},0,X^{n-1},0,X^{n-2},..,0,X^{2},0,X^{1})  \notag \\
&=&\widetilde{I}\ |_{N}(X^{n},0,X^{n-1},0,X^{n-2},..,0,X^{2},0,X^{1})  \notag
\\
&=&I(X^{n},X^{n-1},X^{n-2},..,X^{2},X^{1})
\end{eqnarray}

It is clear that $(ii)$ implies $(i).$

\textbf{Corollary 5.1}. \textit{Let }$I$\textit{\ be a subset of the space }$%
L^{1}(N)$\textit{\ and }$\widetilde{I}$\textit{\ its image by the mapping }$%
\thicksim $\textit{\ such that }$J=\widetilde{I}|_{M}$\textit{\ is an ideal
in }$L^{1}(N),$\textit{\ then the following conditions are verified.}

$(i)$\textit{\ }$J$\textit{\ is a closed ideal in the algebra }$L^{1}(M)$ 
\textit{\ if and only if }$I$\textit{\ is a left closed ideal in the algebra 
}$L^{1}(N).$

$(ii)$\textit{\ }$J$\textit{\ is a maximal ideal in the algebra }$L^{1}(M)$ 
\textit{\ if and only if }$I$\textit{\ is a left maximal ideal in the
algebra }$L^{1}(N).$

$(iii)$\textit{\ }$J$\textit{\ is a prime ideal in the algebra }$L^{1}(M)$ 
\textit{\ if and only if }$I$\textit{\ is a left prime ideal in the algebra }%
$L^{1}(N).$

$(iv)$\textit{\ }$J$\textit{\ is a dense ideal in the algebra }$L^{1}(M)$ 
\textit{\ if and only if }$I$\textit{\ is a left dense ideal in the algebra }
$L^{1}(N).$

The proof of this corollary results immediately from theorem $\mathbf{5.1.}$

\textbf{5.2. Some New Groups}

\textbf{I- }Let $\mathbb{R}_{-}^{\star }=\left\{ x\in \mathbb{R};\text{ }x%
\text{ }\langle 0\text{ }\right\} ,$ with law \ 
\begin{equation}
x\cdot y=-x\text{ }.y
\end{equation}%
for all $x\in $ $\mathbb{R}_{-}^{\star }$ and $y\in \mathbb{R}_{-}^{\star },$
where $\cdot $ signifies the product in $\mathbb{R}_{-}^{\star }$ and $.$
signifies the ordinary product of two real numbers.

($\mathbb{R}_{-}^{\star },$ $\cdot )$\textit{\ with this new law} $\cdot $ 
\textit{becomes a commutative group isomorphic with the multiplicative group 
}$(\mathbb{R}_{+}^{\star },$ $.)$ \textit{and with the real vector group }($%
\mathbb{R}$, $+).$

\textit{Proof: }the identity element is $-1$ because

\begin{eqnarray}
(-1)\cdot x &=&-(-1)\cdot x=x,\text{ }and  \notag \\
x\cdot (-1) &=&-x\text{ }.(-1)=x
\end{eqnarray}

\textbf{\ }If $(x$ ,$y)\in \mathbb{R}_{-}^{\star }\times \mathbb{R}%
_{-}^{\star }$, such that $x\cdot y=-1,$ then we get%
\begin{equation}
x\cdot y=-x\text{ }.y=-1
\end{equation}

From this equation we obtain $y=\frac{1}{x}=x^{-1}\in \mathbb{R}_{-}^{\star
},$ which is the inverse of $x.$ Now let $x,$ $y,$ and $z$ $\ $be three
elements belong to $\mathbb{R}_{-}^{\star },$ then we have 
\begin{equation}
(x\cdot y)\cdot z=-(x\cdot y).z=-(-x\text{ }.y).z=x\text{ }.y.z
\end{equation}%
and%
\begin{equation}
x\cdot (y\cdot z)=x\cdot (-y.z)=-x.(-y.z)=x.y.z
\end{equation}

So the law is associative.

Clearly the law $\cdot $ is commutative because%
\begin{equation}
x\cdot y=-x\text{ }.y=-y.x=y\cdot x
\end{equation}

Let $\psi :$ $\mathbb{R}_{-}^{\star }\rightarrow \mathbb{R}_{+}^{\star }$
the mapping defined by 
\begin{equation}
\psi (x)=-x\text{ }
\end{equation}%
then we get 
\begin{eqnarray}
\psi (x\cdot y) &=&-(x\text{ }\cdot y)  \notag \\
&=&x.y=(-x).(-y)=\psi (x).\psi (y)
\end{eqnarray}
\ 

Another hand we have 
\begin{equation}
\Psi (x)=\Psi (y)\Longrightarrow -x=-y\Longrightarrow x=y
\end{equation}%
and if $y\in \mathbb{R}_{+}^{\star },$ there is $x\in \mathbb{R}_{-}^{\star
},$ such that $\Psi (x)=y.$ In fact it is enough to take $x=-y$

That means $\psi $ is a group isomorphism from $\mathbb{R}_{-}^{\star }$\
onto $\mathbb{R}_{+}^{\star }.$Since the group $(\mathbb{R}_{+}^{\star
},\cdot )$ is isomorphic with the real vector group $(\mathbb{R},+)$ , so
the new group $\mathbb{R}_{-}^{\star }$ is isomorphic with the group $(%
\mathbb{R},+)$

Then $\mathbb{R}^{\star }=$ $\mathbb{R}_{-}^{\star }$\ $\cup $ $\mathbb{R}%
_{+}^{\star }=\mathbb{R}_{+}^{\star }\ \cup \mathbb{R}_{+}^{\star }$

\textbf{II-} \textit{The group} $\mathbb{R}_{+}^{\star }\times \mathbb{R}%
_{-}^{\star }$ , \textit{which is the direct product of} $\ \mathbb{R}%
_{+}^{\star }$ \textit{and} $\mathbb{R}_{-}^{\star }$ \textit{is topological
and algebraic isomorphic onto} $(\mathbb{R}^{2},+)$

\textit{Proof: Let }$\Psi $: $\mathbb{R}_{+}^{\star }\times \mathbb{R}%
_{-}^{\star }\rightarrow \mathbb{R}^{2}$ the mapping defined by%
\begin{equation}
\Psi (x,y)=(\ln x,\ln \left\vert y\right\vert )
\end{equation}

The mapping $\Psi $ is one-to-one and onto. Now, let us to prove that $\Psi $
is group homomorphism

\begin{eqnarray*}
&&\Psi ((x,y)(s,t))=\Psi ((x.s,y\bullet t)=(\ln (x.s),\ln \left\vert
y\bullet t\right\vert ) \\
&=&(\ln (x.s),\ln \left\vert (-1)y.t\right\vert )=(\ln x+\ln s,\ln
\left\vert y.t\right\vert ) \\
&=&(\ln x+\ln s,\ln \left\vert y\right\vert +\ln \left\vert t\right\vert
)=(\ln x,\ln \left\vert y\right\vert )+(\ln s,\ln \left\vert t\right\vert )
\\
&=&\Psi ((x,y)+\Psi (s,t))
\end{eqnarray*}

The inverse of $\Psi $ is given by%
\begin{equation}
\Psi ^{-1}(a,b)=(e^{a},-e^{b})
\end{equation}

\textbf{III- }\textit{The group}\textbf{\ }$\mathbb{R}_{+}^{\star }\times 
\mathbb{R}_{-}^{\star }$ \textit{is also isomorphic onto the group} $(%
\mathbb{C}
,+$\textbf{\ }$)$

\textit{Proof: }It is enough to consider the mapping $\Phi :\mathbb{R}%
_{+}^{\star }\times \mathbb{R}_{-}^{\star }\longrightarrow 
\mathbb{C}
$ defined by $\Phi (x,y)=\ln x+\ln \left\vert y\right\vert ^{i},$ \ $i=(-1)^{%
\frac{1}{2}}$.

\textbf{Acknowledgment. }For all the results of this paper\textbf{\ }the
authors advise the readers to refer my paper $[9],$

\end{document}